\documentclass[twocolumn]{autart}    % Enable this line and disable the 
\usepackage{graphics} % for pdf, bitmapped graphics files
\usepackage{epsfig} % for postscript graphics files
\usepackage{mathptmx} % assumes new font selection scheme installed
\usepackage{times} % assumes new font selection scheme installed
\usepackage{amsmath} % assumes amsmath package installed
\usepackage{amssymb}  % assumes amsmath package installed
\usepackage{xcolor}
\usepackage{subfigure}
\usepackage{tikz}

 \usepackage{natbib}

\allowdisplaybreaks

\newtheorem{theorem}{Theorem}

\newtheorem{lemma}{Lemma}

\begin{document}

\begin{frontmatter}
%\runtitle{Insert a suggested running title}  % Running title for regular 
                                              % papers but only if the title  
                                              % is over 5 words. Running title 
                                              % is not shown in output.

\title{Pitchfork-bifurication-based competitive and collaborative control of an E-bike system\thanksref{footnoteinfo}} % Title, preferably not more 
                                                % than 10 words.

\thanks[footnoteinfo]{An early version of the results presented here were reported, without proofs and experimental evidences of competitive/cooperative scenarios, at the 2019 IEEE 58th Conference on Decision and Control \citep{9029911}. Sweeney and Lhachemi are joint first authors. Russo and Shorten are joint last authors. }

\author[Shaun]{Shaun Sweeney}\ead{s.sweeney21@imperial.ac.uk},    % Add the 
\author[Hugo]{Hugo Lhachemi}\ead{hugo.lhachemi@centralesupelec.fr},               % e-mail address 
\author[Andrew]{Andrew Mannion}\ead{andrew.mannion@ucdconnect.ie},
\author[Giovanni]{Giovanni Russo}\ead{giovarusso@unisa.it},  % (ead) as shown
\author[Bob]{Robert Shorten}\ead{r.shorten@imperial.ac.uk}

\address[Shaun]{Department of Electrical and Electronic Engineering, Imperial College London, UK}  % Please supply                                              
\address[Hugo]{Universit{\'e} Paris-Saclay, CNRS, CentraleSup{\'e}lec, Laboratoire des signaux et syst{\`e}mes, 91190, Gif-sur-Yvette, France}             % full addresses
\address[Andrew]{Department of Electronic and Electrical Engineering, University College Dublin, Ireland.}  
\address[Giovanni]{Department of Information and Electrical Engineering and Applied Mathematics, University of Salerno, Italy}        % here.
\address[Bob]{Dyson School of Design Engineering, Imperial College London, UK}        % here.

                             % keyword list or with the 
                                          % help of the Automatica 
                                          % keyword wizard

\begin{abstract}                          % Abstract of not more than 200 words.
This paper is concerned with the design of a human-in-the-loop system for deployment on a smart pedelec (e-bike). From the control-theoretic perspective, the goal is not only to use the power assistance of the e-bike to reject disturbances along the route but also to manage the possibly competitive interactions between a human and the motor intervention. Managing the competitive/cooperative nature of the interactions is crucial for applications in which we wish to control physical aspects of the cycling behavior (e.g. heart rate and breathing rate).  The basis of the control is a pitchfork bifurcation system, modeling the interactions,  augmented using ideas from gain-scheduling. In vivo experiments have been conducted, showing the effectiveness of the proposed control strategy. 
\end{abstract}

\end{frontmatter}

\section{Introductory remarks}

Over the last few years, electric bikes have become an increasingly popular mode of transportation across many jurisdictions~\citep{booke}. Key reasons for this trend include the fact these {\em e-bikes} are environmentally friendly, enable the design of shared transportation services \citep{Sharing_book}, thus alleviating the  stress on more conventional transportation modes caused by the enormous challenges currently faced by urban areas \citep{Seto16083}. Indeed, the benefits of convincing citizens to transition from cars to e-bikes, include reduced noise and particle emissions from road vehicles, as well as reduced road and parking congestion.  Moreover,  the electric motor on the e-bikes can assist cyclists in completing their journeys, thus helping to overcome the traditional barriers to active transportation of age, topology and wind. Indeed, an impact of the Covid-19 pandemic is that the e-bike revolution is likely to further accelerate as consumers shy away from modes of public transportation due to virus-angst. A further advantage of e-bikes over competitors such as e-scooters is not only the active component of the travel, but also that they respect the separation between pedestrian and road users; simply put, cycling is not ever tolerated on footpaths, whereas social conventions appear to tolerate the use of high power scooters in pedestrian zones.\newline

E-bikes are very attractive, not only due to the electrical assist that can be used to make cycling easier for the cyclist, but also because the presence of a battery offers a host of opportunities to develop services that can be used for benefit of cyclists. and even other road users. Examples of such services include collaborative sensing applications to make cycling safer for cyclists, applications to promote circularity in the use of the e-bike, triage applications to assist emergency services, and applications to protect other vulnerable road users (such as pedestrians). 
From a control perspective, the electric motor offers a really exciting actuation possibility. Not only can be this be used to reject disturbances (wind, hills), but also to provide new services to the cyclist. These  range from basic routing services to avoid e.g. polluted areas and minimize travel time \citep{sweeney2018context,8317888}, but the battery can be used to provide more advanced actuation services such as nudges designed to manage {\em how} the e-bike power-assist can {\em modify} the cyclist behavior during his/her trip. These latter applications classes are challenging for a number of reasons.  Consider, for example,  the problem of designing a control algorithm able to regulate the ventilation rate of the  cyclist in areas of elevated pollution \citep{sweeney2018context}.  To achieve this, the algorithm will need to  judiciously orchestrate  electric assistance.  However, since this reduces the load on the cyclist,  their natural instinct may be to increase cycling effort, thereby possibly increasing his/her ventilation rate and hence missing the control goal.  Thus,  in order to avoid this situation, the control must be designed to account for situations where the cyclist can compete with the motor. 

Hence, a key control-theoretical  challenge is that the actuation provided by the e-bike not only needs to handle disturbances associated to e.g. the topology of the route but it also needs to handle the interaction with the cyclist, which can be both competitive and cooperative.  Motivated by this, we present here a human-in-the-loop control system for e-bikes that is able to manage the interactions between a human and a motor intervention, for applications in which we wish to control physical aspects of the cycling behavior.  In particular, our objective is to develop a gain-scheduling cyber-physical control system~\citep{leith2000survey,rugh2000research} that manages such an interaction between the cyclist and the electrical motor.   We now briefly survey some related works.

%In this latter context we are currently working on a number of applications underpinned by control and dynamical systems theory. These include:
%\begin{itemize}
%\item control strategies to nudge the cyclist to obey traffic light signalling;
%
%\item camera based control strategies that help cyclists avoid pedestrians on cycle paths;
%
%\item control strategies to manage cyclist perspiration so that they do not arrive at work sweating;
%
%\item and control strategies to regulate the breathing rate of the cyclist (to, for example, mitigate the effects of pollution peaks). 
%\end{itemize}

\subsection*{Related work}

Existing literature on e-bikes is mainly devoted to their energy management. Interestingly, as noted in \citep{guanetti2017optimal} most hybrid electric bicycles  on the market determine the level of assistance without explicitly accounting for energy efficiency and/or human behavior. Recent works such as \citep{6859373b,6580849}, presented an offline approach to compute optimal pacing by taking into account muscolar fatigue. We also recall here \citep{6669315,6580369,6314871} which are focused on the design of an online energy management strategy for full hybrid electric bicycles (also experimentally validated).  More recently,  a control based, energy management system for series hybrid electric bicycles is presented in~\citep{guanetti2017optimal}, although without consideration of the competitive interaction between cyclist and the electrical motor system. In a broader context, the development of e-bike services are presented in a number of publications; see~\citep{smaldone,kiefer}, as well as the the work carried out as part of the Copenhagen-Wheel project at MIT, and the recent paper from IBM Research~\citep{IBM2018}.  Other services that can be enabled by e-bikes leverage the fact that these can be used for crowdsourcing \citep{9244209} applications (for example pollution and/or routes monitoring), and to provide other services to municipalities, pedestrians, and road users, in exchange for monetary reward (for example, mobile wifi hubs); see~\citep{parked}. Further, in many of contemporary e-bikes, the battery can be large - up to 500 Whrs. Thus, when aggregated together, groups of batteries can provide energy buffering services to both infrastructure and citizens alike~\citep{dlt_paper}.  Our work is also related to the rich body of work on the design of cyber-physical systems, that is the subject of much interest in the control engineering community. For context, we refer the interested reader to the the recently published paper on human-in-the-loop control systems ~\citep{inoue2019weak} and the excellent recent review paper~\citep{hirsch}. We note briefly that while conceptually the work presented here is related to this latter body of work, the {\em competitive}/{\em cooperative} nature of the control is a somewhat unique feature of the problem class considered in this paper. Perhaps, as such, our work is most related to the stabilisation of interacting unstable systems discussed by Narendra~\citep{Narendra} and the agent coexistence problems discussed in the networking community~\citep{rade}; our solution owes much inspiration to these latter two works. We also note that the work presented in this paper builds on our work presented in~\citep{sweeney2018context}. In this work an e-bike design is presented, as well as a holistic control strategy (without proofs) managing the interaction between the cyclist and the e-bike motor. 

\subsection*{Statement of contributions}

In this paper we present a human-in-the-loop control system for an e-bike that is able to manage the competitive/cooperative nature of the interactions between a cyclist and the e-bike. To the best of our knowledge, this is the first control system for e-bikes that explicitly takes into account the possibly non-cooperative nature of humans when provided with power assistance from the e-bike.  We present a control system based on the use of gain-scheduling and prove a number of closed-loop properties.  While the work presented in this paper builds on~\citep{sweeney2018context}, in the latter work we only presented a principled control design, without giving any guarantees on its closed-loop performance.  In the context of the above literature, our contributions can be summarized as follows:
\begin{itemize}
\item we present a novel control algorithm to manage the interactions between the cyclist and the e-bike. The key idea behind the algorithm, which is also a technical contribution of this paper, is based on the use a full non-linear gain scheduling~\citep{leith2000survey,rugh2000research}  design based on a pitchfork bifurcation.  As we shall see,  this leads to the design a control system that can handle a full envelope of operating conditions;
\item we prove certain properties of the closed-loop system. In particular, we give a rigorous assessment of the robustness and tracking performance (in a sense specified in Section \ref{sec: control design});
\item we investigate the effectiveness of the closed-loop system via {\em in-vivo} experiments. Namely, with the experiments, we test the design on a real e-bike in a controlled set-up.
\end{itemize}

\section{System description}
\begin{figure}[h]
\begin{center}
{\includegraphics[width=3in]{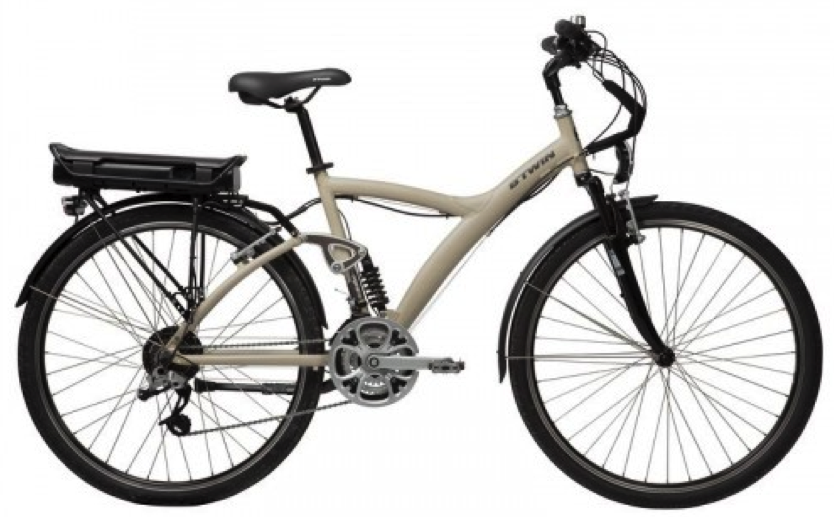}}
\caption{Electric bike from BTwin.}
\label{Ebike}
\end{center}
\end{figure}

We use the infrastructure from~\citep{sweeney2018context}, which is described here for completeness.  The e-bike that we use is a modified BTwin Original 700 purchased from Decathlon (see Figure \ref{Ebike}).  While the original e-bike was equipped with a Samsung Li-ion 36V battery and a Bafang controller, we modified it in several ways:
\begin{itemize}
\item in order to facilitate the actual design and implementation of the control algorithm,  we replaced the original motor controller with the Grinfineon C4820-GR.  This appeared to be a more advanced controller hardware, suitable for our purposes;
\item several measurement sensors were added to the e-bike. These additional sensors were aimed at measuring pedal torque/speed using a THUN X-CELL RT sensor\footnote{http://www.ebikes.ca/shop/electric-bicycle-parts/torque-sensors/thun-120l.html}), battery voltage, motor current, wheel speed, motor temperature. We also included in the design  brake and hand throttle sensors.
\end{itemize}

We also developed the hardware and software infrustructure to process the data from sensors. In particular, all the sensors were read by either an Arduino (brake and hand throttle sensors) or by an off-the-shelf computer system. In the experiments, we used the Cycle Analyst\protect\footnote{www.ebikes.ca} to retrieve the data from the sensors. In both cases, data would be communicated to a smartphone app (that we also designed) using a 
bespoke especially designed Arduino-controlled Bluetooth module. Control inputs were  sent to the e-bike controller using the same Bluetooth based communication system. 

Moreover,  in order to guarantee {\em reproducibility}  of our experimental set-up, we decided to perform tests in a controlled environment that would allow to reproduce all the operating conditions described in Section \ref{sec: results}.  

\section{Control design}\label{sec: control design}

In this section, we describe the control strategy adopted in this work. First, we present a simple model of the e-bike power-generation process. Second, we discuss the necessity of a human-in-the-loop-based design of the control strategy for fostering a cooperative behavior of the cyclist. Third, the basic idea used to realise such a control design  is presented. Finally, the fully functional gain-scheduling control design procedure is introduced and theoretical analysis of the design  presented.

In the sequel, $\mathbb{R}_+$ and $\mathbb{R}_+^*$ denote the sets of non-negative real and positive real numbers, respectively. For any given connected set $S \subset \mathbb{R}^n$ and any interval $I \subset \mathbb{R}$, $\mathcal{C}^0(S;I)$ denotes the set of functions $g : S \rightarrow I$ that are continuous. Similarly,  $\mathcal{C}^1(S;\mathbb{R})$ denotes the set of functions $g : S \rightarrow \mathbb{R}$ that are continuously differentiable.

\subsection{Basic modeling}

As in \citep{sweeney2018context}  we make the two following assumptions. First, we can neglect any dynamics associated with the e-bike's electric motor (motor in short). Second, we consider a fixed mechanical gear setting.  The first instantaneous input power to be considered is from the motor to the bike $P_{M_\mathrm{in}}(t)$. The control input to the motor is denoted by $Y(t)$ and is related to the motor current $I_M(t)$ via $I_M(t) = \mu Y(t)$ where $\mu > 0$ is a constant. Introducing $V_M$ the motor voltage, the power delivered by the motor is given by $P_{M_\mathrm{in}}(t) = V_M I_M(t)$. The second instantaneous input power to be considered is from the human to the bike $P_{H_\mathrm{in}}(t)$. It is expressed as $P_{H_\mathrm{in}}(t) = \tau_p(t) \omega_p(t)$ where $\tau_p(t)$ represents the torque provided by the cyclist at the pedal and $\omega_p(t)$ is the angular velocity at pedals. Introducing $E_m \in (0,1)$ and $E_c \in (0,1)$ the efficiency of the motor and the crankset, respectively, the output motor power and output human power are given by $P_{M_\mathrm{out}}(t) = E_m P_{M_\mathrm{in}}(t)$ and $P_{H_\mathrm{out}}(t) = E_c P_{H_\mathrm{in}}(t)$, respectively. Thus, the total power available to move the rear wheel can be approximated as $P_w(t) = P_{M_\mathrm{out}}(t) + P_{H_\mathrm{out}}(t)$.

\subsection{Control Algorithm}
The control objective is to regulate the fraction of effort delivered by the human. For example, such a control strategy is important to control the transpiration or the ventilation rate of the cyclist. To achieve such a control objective, we introduce the proportion $m \in [0,1]$ of power provided by the cyclist:
\begin{equation}\label{eq: m - naive approach}
m(t) = \dfrac{P_{H_\mathrm{out}}(t)}{P_{M_\mathrm{out}}(t)+P_{H_\mathrm{out}}(t)} .
\end{equation}
A value of $m(t)$ close to $0$ indicates that the system is essentially in electric mode while a value of $m(t)$ close to $1$ means that most of the power is delivered by the human. For a given desired value $m^* \in [0,1]$, the control objective is that $m$ achieves the setpoint tracking of $m^*$. \newline

Note that the configuration $m = 0$ with $P_{H_\mathrm{out}} = 0$ and $P_{M_\mathrm{out}} > 0$ corresponds to a full electrical mode and is undesirable as it would transform the e-bike into a motorbike, which is illegal in some jurisdictions. Therefore, we assume the existence of a constant $\eta \in (0,1)$ such that $m^*(t) \in [\eta , 1]$ for all $t \geq 0$.\newline

In practice, the direct use of the quantities $P_{M_\mathrm{out}}(t)$ and $P_{H_\mathrm{out}}(t)$ is not adequate because of biases, uncertainties in the system operation, and the stochasticity of the cyclist behavior over short period of time (such as sudden speed variations over short distances). Consequently, the direct control of $m$ defined by (\ref{eq: m - naive approach}) can be counterproductive as it might lead to a very erratic control effort $Y$. This would result in a very erratic power delivery of the motor $P_{M_\mathrm{out}}$, inducing safety issues and a lack of comfort for the cyclist. Consequently, we introduce the following filtered version of the proportion $m \in [0,1]$ of power provided by the cyclist:
\begin{equation}
m(t) = \dfrac{\overline{P}_{H_\mathrm{out}}(t)}{\overline{P}_{M_\mathrm{out}}(t)+\overline{P}_{H_\mathrm{out}}(t)} ,
\end{equation}
where $\overline{P}_{H_\mathrm{out}}$ and $\overline{P}_{M_\mathrm{out}}$ are the corresponding filtered and averaged versions of $P_{H_\mathrm{out}}$ and $P_{M_\mathrm{out}}$ (see~\citep{shaun2017} for details). In this setting, the control objective remains that $m$ should track a given reference value $m^*$. For notational simplicity, $\overline{P}_{H_\mathrm{out}}$ and $\overline{P}_{M_\mathrm{out}}$ are simply denoted by $P_{H_\mathrm{out}}$ and $P_{M_\mathrm{out}}$ in the sequel. \newline

Denoting by $\Delta T$ the sampling time and introducing the regulation error $e_k = m^*(k \Delta T) - m(k \Delta T)$, it was proposed and implemented in~\citep{sweeney2018context} the following discrete time (with zero order holder) integral controller:
\begin{equation}
Y_{k+1} = Y_{k} - \gamma e_k ,
\end{equation}
where $Y_{k} = Y(k \Delta T)$ and $\gamma > 0$ is a proportional gain. The validity of this approach was assessed in~\citep{sweeney2018context} by means of experiments. However, it does not fully account for the human-in-the-loop aspect of the problem, i.e., the complex interactions between the cyclist and the motor. Specifically, a decrease of the value of the reference $m^*$ should reduce the workload of the cyclist by decreasing accordingly the value of $m$. However, due to this reduction of the workload on the cyclist, their instinct may be to increase cycling effort, which is the opposite of the original objective. Consequently, the control strategy must consider the possibility of such an uncooperative behavior of the cyclist by progressively switching-off the motor.

\subsection{Basic idea}
\begin{figure}[htb]
\begin{tikzpicture}
\draw[->] (-0.5,0) -- (6.25,0) node[below] {$P_{H_\mathrm{out}}$};
\draw[->] (0,-0.5) -- (0,3) node[left] {$\sqrt{f(P_{H_\mathrm{out}})}$};
\draw[scale=1,domain=0:3,smooth,variable=\x,blue] plot ({\x},{0.8*\x});
\draw[scale=1,domain=3:6,smooth,variable=\x,red] plot ({\x},{2.4*exp(1.5*(3-\x))});
\draw[dashed] (3,0) node[below] {$P_{T}$} -- (3,2.85) ;
\node[align=left,text width=2.5cm] at (1.5,2.5) {\small Human and motor \textcolor{blue}{cooperate}};
\node[align=left,text width=2.5cm] at (4.75,2.5) {\small Human and motor \textcolor{red}{compete}};
\end{tikzpicture}
\caption{Typical function $f$ specifying cooperative and competitive behaviors}
\label{fig: function f coop vs comp}
\end{figure}
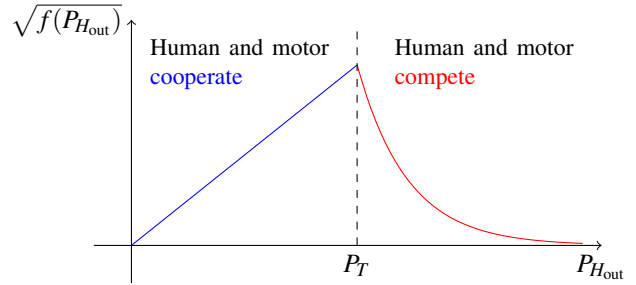

As suggested in~\citep{sweeney2018context}, one approach for tackling this problem consists in the following control strategy:
\begin{equation}\label{eq: pitchfork bifurcation - basic idea}
\dot{P}_{M_\mathrm{out}} = f \left( P_{H_\mathrm{out}} \right) P_{M_\mathrm{out}} - P_{M_\mathrm{out}}^3 ,
\end{equation}
which is the normal form of a supercritical pitchfork bifurcation~\citep{kuznetsov2013elements}. The function $f$ is designed based on engineering considerations such that $f(0) = 0$ and $f(x) > 0$ for any $x > 0$. The motivation behind the use of (\ref{eq: pitchfork bifurcation - basic idea}) is to regulate the interaction between the cyclist and the motor. 
\begin{itemize}
\item When $P_{H_\mathrm{out}} = 0$, we have $f(P_{H_\mathrm{out}}) = 0$. Hence the only equilibrium point of (\ref{eq: pitchfork bifurcation - basic idea}) is $P_{M_\mathrm{out}} = 0$, which in turn is globally asymptotically stable. That is $P_{M_\mathrm{out}} \rightarrow 0$ when $t \rightarrow + \infty$, whatever the initial condition $P_{M_\mathrm{out}}(0) \geq 0$ is.
\item Whenever\footnote{Subscript "e" stands for an equilibrium value.} $ P_{H_\mathrm{out}}(t) =  P_{H_\mathrm{out},e} > 0$, the origin of (\ref{eq: pitchfork bifurcation - basic idea}) becomes unstable and the new stable equilibrium $P_{M_\mathrm{out},e} = \sqrt{f(P_{H_\mathrm{out},e})}$ appears. In this case, $P_{M_\mathrm{out}}$ converges to $P_{M_\mathrm{out},e}$ as $t \rightarrow + \infty$, whatever the initial condition $P_{M_\mathrm{out}}(0) > 0$ is.
\end{itemize}
In this context, the shape of the function $f$ is tuned such that the motor assists the cyclist. Specifically, $m$ tracks $m^*$ as long as $P_{H_\mathrm{out}}$ is below a given threshold $P_T > 0$, while the motor is gradually switched off when $P_{H_\mathrm{out}} > P_T$. The graph of such a typical function $f$ is depicted in Fig.~\ref{fig: function f coop vs comp}.

\subsection{Fine tuning of the pitchfork bifurcation by means of gain-scheduling control design}

The control law adopted in this paper takes the following form:
\begin{subequations}
\begin{align}
\dot{P}_{M_\mathrm{out}} & = \alpha_{m^*} \left( P_{H_\mathrm{out}} \right) \left[ f_{m^*} \left( P_{H_\mathrm{out}} \right) P_{M_\mathrm{out}} - P_{M_\mathrm{out}}^3 \right] \label{eq: pitchfork bifurcation} \\
P_{M_\mathrm{out}}(0) & = P_{M_{\mathrm{out}},0} \geq 0 \label{eq: pitchfork bifurcation - IC}
\end{align}
\end{subequations}
with $ \alpha_{m^*} \left( P_{H_\mathrm{out}} \right) > 0$, $f_{m^*} ( 0 ) = 0$, and $f_{m^*} \left( P_{H_\mathrm{out}} \right) > 0$ for all $P_{H_\mathrm{out}} > 0$. Functions $(P_{H_\mathrm{out}},m^*) \mapsto f_{m^*}\left( P_{H_\mathrm{out}} \right)$ and $(P_{H_\mathrm{out}},m^*) \mapsto \alpha_{m^*}\left( P_{H_\mathrm{out}} \right)$, that are selected in $\mathcal{C}^0(\mathbb{R}_+ \times [\eta,1] ; \mathbb{R}_+)$, must be tuned in function of $m^* \in [\eta , 1]$. The former aims at capturing the desired trade-off between cooperative and competitive behaviors of the cyclist, while the latter is for the tuning of the speed of convergence of the control strategy. \newline

The function $(P_{H_\mathrm{out}},m^*) \mapsto f_{m^*}\left( P_{H_\mathrm{out}} \right)$ is tuned such that, for any constant reference input $m^*$ and any equilibrium value $P_{H_\mathrm{out},e} \leq P_T(m^*)$, we have at the equilibrium that
\begin{equation*}
m^* 
= m_e 
= \dfrac{P_{H_\mathrm{out},e}}{P_{M_\mathrm{out},e} + P_{H_\mathrm{out},e}}
= \dfrac{P_{H_\mathrm{out},e}}{\sqrt{f_{m^*}(P_{H_\mathrm{out},e})} + P_{H_\mathrm{out},e}} .
\end{equation*}
In this work, this is achieved by defining for $P_{H_\mathrm{out}} \leq P_T(m^*)$:
\begin{equation}\label{eq: def fm*}
f_{m^*}\left( P_{H_\mathrm{out}} \right) 
=
\left[ \dfrac{1 - m^*}{m^*} P_{H_\mathrm{out}} \right]^2 
\end{equation}
and for $P_{H_\mathrm{out}} > P_T(m^*)$:
\begin{align*}
f_{m^*}\left( P_{H_\mathrm{out}} \right)
& =
\left[ \dfrac{1 - m^*}{m^*} P_T(m^*) \right]^2 \\
& \phantom{=}\; \times \left[ 1 + 2 \left( \gamma + \dfrac{1}{P_T(m^*)} \right) ( P_{H_\mathrm{out}} - P_T(m^*) ) \right] \\
& \phantom{=}\; \times e^{-2\gamma \left[ P_{H_\mathrm{out}} - P_T(m^*) \right]} ,
\end{align*}
where $\gamma >0$ is the exponential decay rate for the competitive scenario. $P_T(m^*) > 0$ denotes the threshold associated with the current value of the reference signal $m^*\in [\eta , 1]$. In other words, for a given value of $m^*$, $P_T(m^*) > 0$ stands for the maximum value of the effort provided by the cyclist that is considered as a cooperative behavior. For large values of $m^*$ (i.e., close to $1$), the cyclist should be able to modulate freely his effort. Therefore, the corresponding value of the threshold $P_T(m^*)$ should be large enough. Conversely, the effort of the cyclist should be constrained when $m^*$ is close to the lower bound $\eta$. So, in this case, the value of the threshold $P_T(m^*)$ should be reduced. Thus, we select $m^* \mapsto P_T(m^*)$ as an element of $\mathcal{C}^0([\eta,1] ; \mathbb{R}_+^*) \cap \mathcal{C}^1([\eta,1] ; \mathbb{R})$ satisfying $\dot{P}_T > 0$, i.e., is a (strictly) positive and increasing function of the reference value $m^* \in [\eta , 1]$. Thus, it is easy to see that $(P_{H_\mathrm{out}},m^*) \mapsto f_{m^*}\left( P_{H_\mathrm{out}} \right)$ is an element of $\mathcal{C}^0(\mathbb{R}_+ \times [\eta,1] ; \mathbb{R}_+) \cap \mathcal{C}^1(\mathbb{R}_+ \times [\eta,1] ; \mathbb{R})$.  \newline

We now tune the function $(P_{H_\mathrm{out}},m^*) \mapsto \alpha_{m^*}\left( P_{H_\mathrm{out}} \right)$ to tune the convergence rate of the control law. To do so, assuming a constant reference signal $m^*$, we resort to the following Jacobian linearization of the controller dynamics (\ref{eq: pitchfork bifurcation}) in the vicinity of the equilibrium point $P_{M_\mathrm{out},e} = \sqrt{f_{m^*}(P_{H_\mathrm{out},e})}$:
\begin{equation*}
\delta\dot{P}_{M_\mathrm{out}} = - 2 \alpha_{m^*}(P_{H_\mathrm{out},e}) f_{m^*}(P_{H_\mathrm{out},e}) \delta P_{M_\mathrm{out}} .
\end{equation*}
In order to impose a constant convergence rate $\kappa > 0$ in the vicinity of any equilibrium point, we select $\alpha_{m^*}\left( P_{H_\mathrm{out}} \right) = \kappa / ( 2 f_{m^*}(P_{H_\mathrm{out}}) ) > 0$ which yields $\delta\dot{P}_{M_\mathrm{out}} = - \kappa \delta P_{M_\mathrm{out}}$. However, for either small or large values of $P_{H_\mathrm{out}}$, $f_{m^*}(P_{H_\mathrm{out}})$ is getting arbitrarily close to zero and thus induces the divergence of $\alpha_{m^*}\left( P_{H_\mathrm{out}} \right)$. To avoid this pitfall, we introduce a threshold $\varepsilon_T > 0$ and set $\alpha_{m^*} = \kappa / ( 2 \epsilon_T )$ for any $P_{H_\mathrm{out}}$ satisfying $f_{m^*}(P_{H_\mathrm{out}}) \leq \epsilon_T$. Thus, we obtain that $\alpha_{m^*}\left( P_{H_\mathrm{out}} \right) = \kappa / ( 2 \max( f_{m^*}(P_{H_\mathrm{out}}) , \epsilon_T ) ) > 0$ which is bounded above by the constant $\kappa / ( 2 \epsilon_T )$.

\subsection{Robustness of the control strategy}

In order to validate the proposed controlled strategy, we derive in this section a tracking performance of $m$ with respect to the reference signal $m^*$.

\subsubsection{Preliminary results}

We introduce\footnote{It can be easily seen from the definition of $f_{m^*}(P_{H_\mathrm{out}})$ that  $C_f < + \infty$ by using the fact that $0 < P_T(\eta) \leq P_T(m^*) \leq P_T(1) < +\infty$ for all $m^* \in [\eta , 1]$ and that $x e^{- a x} \leq e^{-1}/a$ for all $x \geq 0$ and $a > 0$.}
\begin{equation*}
C_f \triangleq \sup \limits_{(m^*,P_{H_\mathrm{out}}) \in [\eta , 1] \times \mathbb{R}_+}\; f_{m^*}(P_{H_\mathrm{out}}) 
\in \mathbb{R}_+^* 
\end{equation*}
and, for any constants $P_{H,M} > P_{H,m} > 0$, we define\footnote{$c_f(P_{H,M},P_{H,m}) > 0$ because it is defined as the infimum over a compact set of a continuous function that is (strictly) positive.}
\begin{equation*}
c_f(P_{H,M},P_{H,m}) \triangleq \inf\limits_{m^* \in [\eta , 1] \,,\, P_{H,m} \leq P_{H_\mathrm{out}} \leq P_{H,M}}\; f_{m^*}(P_{H_\mathrm{out}}) \in \mathbb{R}_+^* .
\end{equation*} 

We can now state the following result dealing with the maximum amplitude of the power delivered by the motor for a suitable range of admissible initial conditions $P_{M_{\mathrm{out}},0}$.\newline

\begin{lemma}\label{lem 1}
Let $P_{M_{\mathrm{out}},0} \in [0,\sqrt{C_f}]$, $P_{H_\mathrm{out}} \in \mathcal{C}^0(\mathbb{R}_+;\mathbb{R}_+)$, and $m^* \in \mathcal{C}^0(\mathbb{R}_+;[\eta , 1])$ be given. The unique maximal solution $P_{M_\mathrm{out}}$ of (\ref{eq: pitchfork bifurcation}-\ref{eq: pitchfork bifurcation - IC}) is defined over $\mathbb{R}_+$ and satisfies $P_{M_\mathrm{out}}(t) \in [0,\sqrt{C_f}]$ for all $t \geq 0$.
\end{lemma}

\begin{pf*}{Proof.}
As $t \mapsto \alpha_{m^*(t)} \left( P_{H_\mathrm{out}}(t) \right)$ and $t \mapsto f_{m^*(t)} \left( P_{H_\mathrm{out}}(t) \right)$ are bounded, it is easy to see that the vector field associated with (\ref{eq: pitchfork bifurcation}) is locally lipschitz in $P_{M_\mathrm{out}}$, uniformly in $t \geq 0$. Thus, the Cauchy-Lipschitz theorem~\citep{khalil2002nonlinear} ensures the existence of a unique maximal solution $P_{M_\mathrm{out}}$, defined over $[0,t_M)$ with $0 < t_M \leq + \infty$, associated with the initial condition $P_{M_{\mathrm{out}},0}$.

Noting that $P_{M_\mathrm{out}} = 0$ is the unique solution of (\ref{eq: pitchfork bifurcation}-\ref{eq: pitchfork bifurcation - IC}) associated with the zero initial condition, we immediately deduce that $P_{M_\mathrm{out}}(0) \geq 0$ implies $P_{M_\mathrm{out}}(t) \geq 0$ for all $t \in [0,t_M)$. Let $\epsilon > 0$ be arbitrarily given. We prove by contradiction that $\sup\limits_{t \in [0,t_M)}\; P_{M_\mathrm{out}}(t) \leq \sqrt{C_f + \epsilon}$. Indeed, assuming that $\sup\limits_{t \in [0,t_M)}\; P_{M_\mathrm{out}}(t) > \sqrt{C_f + \epsilon}$, we can introduce $t_0 \triangleq \inf \{ t \in [0,t_M) \,:\, P_{M_\mathrm{out}}(t) > \sqrt{C_f + \epsilon} \} < t_M$. As $P_{M_\mathrm{out}}(0) \leq \sqrt{C_f}$, a continuity argument shows that $t_0 > 0$ and $P_{M_\mathrm{out}}(t_0) = \sqrt{C_f+\epsilon}$. From (\ref{eq: pitchfork bifurcation}), we obtain that 
\begin{align*}
& \dot{P}_{M_\mathrm{out}}(t_0) \\
& = \alpha_{m^*(t_0)} \left( P_{H_\mathrm{out}}(t_0) \right) \left[ f_{m^*(t_0)} \left( P_{H_\mathrm{out}}(t_0) \right) P_{M_\mathrm{out}}(t_0) - P_{M_\mathrm{out}}(t_0)^3 \right] \\
& = \alpha_{m^*(t_0)} \left( P_{H_\mathrm{out}}(t_0) \right) \sqrt{C_f + \epsilon} \left\{ f_{m^*(t_0)} \left( P_{H_\mathrm{out}}(t_0) \right) - (C_f + \epsilon) \right\} \\
& < 0 ,
\end{align*}
where the inequality holds true because $\alpha_{m^*(t_0)} \left( P_{H_\mathrm{out}}(t_0) \right) > 0$ and $f_{m^*(t_0)}(P_{H_\mathrm{out}}(t_0)) \leq C_f$. By continuity, there exists $\delta \in ( 0 , \min(t_M - t_0,t_0))$ such that $\dot{P}_{M_\mathrm{out}}(t) < 0$ for $\vert t - t_0 \vert < \delta$. We deduce that $P_{M_\mathrm{out}}(t) \leq P_{M_\mathrm{out}}(t_0) = \sqrt{C_f+\epsilon}$ for $t_0 \leq t < t_0 + \delta$, which is not consistent with the definition of $t_0$. By contradiction, the claimed result $\sup\limits_{t \in [0,t_M)}\; P_{M_\mathrm{out}}(t) \leq \sqrt{C_f + \epsilon} < + \infty$ holds. We deduce that $t_M = +\infty$ and $\sup\limits_{t \geq 0}\; P_{M_\mathrm{out}}(t) \leq \sqrt{C_f + \epsilon}$. By letting $\epsilon \rightarrow 0^+$, the proof is complete.
\end{pf*} 

Now, assuming that a minimal level of power is provided by the human, we obtain the existence of a lower bound on the power delivered by the motor for a suitable range of admissible initial conditions $P_{M_{\mathrm{out}},0}$.  

\begin{lemma}\label{lem 2}
Let $P_{H,M} > P_{H,m} > 0$ be given constants and define $c_f = c_f(P_{H,M},P_{H,m}) > 0$. Let $P_{M_{\mathrm{out}},0} \in [\sqrt{c_f},\sqrt{C_f}]$, $P_{H_\mathrm{out}} \in \mathcal{C}^0(\mathbb{R}_+;\mathbb{R}_+)$ with $P_{H,m} \leq P_{H_\mathrm{out}}(t) \leq P_{H,M}$ for all $t \geq 0$, and $m^* \in \mathcal{C}^0(\mathbb{R}_+;[\eta , 1])$ be given. Then, the associated trajectory $P_{M_\mathrm{out}}$ of (\ref{eq: pitchfork bifurcation}-\ref{eq: pitchfork bifurcation - IC}) satisfies $P_{M_\mathrm{out}}(t) \geq \sqrt{c_f}$ for all $t \geq 0$.
\end{lemma}

\begin{pf*}{Proof.}
From (\ref{eq: pitchfork bifurcation}) and as $P_{M_\mathrm{out}}(t) \geq 0$, we have that 
\begin{align*}
\dot{P}_{M_\mathrm{out}}(t)
& \geq \alpha_{m^*(t)} \left( P_{H_\mathrm{out}}(t) \right) P_{M_\mathrm{out}}(t) \\
& \phantom{\geq} \times \left( \sqrt{c_f} + P_{M_\mathrm{out}}(t) \right) \left( \sqrt{c_f} - P_{M_\mathrm{out}}(t) \right)  .
\end{align*}
Let $\epsilon \in (0 , \sqrt{c_f})$ be arbitrarily given. Assume that there exists $t \geq 0$ such that $P_{M_\mathrm{out}}(t) < \sqrt{c_f} - \epsilon$. So, we can introduce $t_0 \triangleq \inf \{ t \geq 0 \,:\, P_{M_\mathrm{out}}(t) < \sqrt{c_f} - \epsilon \}$. As $P_{M_\mathrm{out}}(0) \geq \sqrt{c_f}$, we obtain by continuity that $t_0 > 0$ and $P_{M_\mathrm{out}}(t_0) = \sqrt{c_f} - \epsilon$. Consequently, we have that
\begin{align*}
\dot{P}_{M_\mathrm{out}}(t_0)
& \geq \epsilon \alpha_{m^*(t_0)} \left( P_{H_\mathrm{out}}(t_0) \right) \left( 2 \sqrt{c_f} - \epsilon \right) (\sqrt{c_f} - \epsilon) > 0 .
\end{align*}
Thus, there exists $\delta > 0$ such that $P_{M_\mathrm{out}}(t) \geq P_{M_\mathrm{out}}(t_0) = \sqrt{c_f} - \epsilon$ for $t_0 \leq t < t_0 + \delta$, which is not consistent with the definition of $t_0$. By contradiction, we obtain that $P_{M_\mathrm{out}}(t) \geq \sqrt{c_f} - \epsilon$ for all $t \geq 0$. By letting $\epsilon \rightarrow 0^+$, the claimed conclusion holds.
\end{pf*}

\subsubsection{Assessment of the robust tracking of the reference signal}

We can now introduce the main theoretical result of this paper. Specifically, assuming that a minimal level of power is provided by the human, we assess for a suitable range of initial conditions and in the case of a cooperative behavior that $m$ tracks $m^*$.

\begin{theorem}
Let $P_{H,M} > P_{H,m} > 0$ and $p \in \{1 , 2 \}$ be given constants and define $c_f = c_f(P_{H,M},P_{H,m})$. There exist constants $\beta,C_1,C_2,C_3,C_4,C_5 > 0$ such that, for any initial condition $P_{M_{\mathrm{out}},0} \in [\sqrt{c_f},\sqrt{C_f}]$, any $P_{H_\mathrm{out}} \in \mathcal{C}^0(\mathbb{R}_+;\mathbb{R}_+) \cap \mathcal{C}^1(\mathbb{R}_+;\mathbb{R})$ with $P_{H,m} \leq P_{H_\mathrm{out}}(t) \leq P_{H,M}$ for all $t \geq 0$, and any $m^* \in \mathcal{C}^0(\mathbb{R}_+;[\eta , 1]) \cap \mathcal{C}^1(\mathbb{R}_+;\mathbb{R})$, the associated trajectory $P_{M_\mathrm{out}}$ of (\ref{eq: pitchfork bifurcation}-\ref{eq: pitchfork bifurcation - IC}) satisfies the estimate: 
\begin{align}
& \left\vert P_{M_\mathrm{out}}(t) - \sqrt{f_{m^*(t)} \left( P_{H_\mathrm{out}}(t) \right)} \right\vert \nonumber \\
& \qquad\leq
\left\vert P_{M_\mathrm{out},0} - \sqrt{f_{m^*(0)} \left( P_{H_\mathrm{out}}(0) \right)} \right\vert e^{-\beta t} \label{eq: ISS P_M_out} \\
& \qquad\phantom{\leq} 
+ C_1 \sup\limits_{\tau \in [0,t]} \vert \dot{m}^*(\tau) \vert^{1/p}
+ C_2 \sup\limits_{\tau \in [0,t]} \left\vert \dot{P}_{H_\mathrm{out}}(\tau) \right\vert^{1/p} \nonumber
\end{align}
with $\sqrt{c_f} \leq P_{M_\mathrm{out}}(t) \leq \sqrt{C_f}$ for all $t \geq 0$. Furthermore, assuming a cooperative behavior of the cyclist, i.e., $P_{H_\mathrm{out}}(t) \leq P_T(m^*(t))$ for all $t \geq 0$, we have that
\begin{align}
\left\vert m(t) - m^*(t) \right\vert
& \leq
C_3 \left\vert m(0) - m^*(0)  \right\vert e^{-\beta t} \label{eq: ISS m} \\
& \phantom{\leq} 
+ C_4 \sup\limits_{\tau \in [0,t]} \vert \dot{m}^*(\tau) \vert^{1/p}
+ C_5 \sup\limits_{\tau \in [0,t]} \left\vert \dot{P}_{H_\mathrm{out}}(\tau) \right\vert^{1/p} . \nonumber
\end{align}
\end{theorem}

\begin{pf*}{Proof.}
We define for all $t \geq 0$,
\begin{equation*}
V(t) = \dfrac{1}{2} \left( P_{M_\mathrm{out}}(t) - \sqrt{f_{m^*(t)} \left( P_{H_\mathrm{out}}(t) \right)} \right)^2 .
\end{equation*}
As $m^*(t) \in [\eta , 1]$ and $P_{H_\mathrm{out}}(t) \in [ P_{H,m} , P_{H,M} ]$, then $f_{m^*(t)} \left( P_{H_\mathrm{out}}(t) \right) > 0$ for all $t \geq 0$ and thus $V$ is continuously differentiable with:
\begin{align}
\dot{V}
& =  - 2 \alpha_{m^*} \left( P_{H_\mathrm{out}} \right) P_{M_\mathrm{out}} \left( P_{M_\mathrm{out}} + \sqrt{f_{m^*} \left( P_{H_\mathrm{out}} \right)} \right) V \label{eq: dot_V} \\
& \phantom{=}\, - \left\{
\dfrac{\partial f_{m^*}(P_{H_\mathrm{out}})}{\partial m^*} \dot{m}^*
+ \dfrac{\partial f_{m^*}(P_{H_\mathrm{out}})}{\partial P_{H_\mathrm{out}}} \dot{P}_{H_\mathrm{out}}
\right\} \nonumber \\
& \phantom{=\,-}\;\; \times \dfrac{P_{M_\mathrm{out}} - \sqrt{f_{m^*}\left( P_{H_\mathrm{out}} \right)}}{2 \sqrt{f_{m^*}\left( P_{H_\mathrm{out}} \right)}} . \nonumber
\end{align}
We introduce the following quantity\footnote{$F_1,F_2 < + \infty$ because they are defined as the supremum over a compact set of a continuous function. $\zeta > 0$ for the same reason as $c_f(P_{H,M},P_{H,m}) > 0$, see previous footnote.}:
\begin{align*}
F_1 & = \sup\limits_{m^* \in [\eta , 1] \,,\, P_{H,m} \leq P_{H_\mathrm{out}} \leq P_{H,M}}\; \left\vert \dfrac{\partial f_{m^*}(P_{H_\mathrm{out}})}{\partial m^*} \right\vert < +\infty , \\
F_2 & = \sup\limits_{m^* \in [\eta , 1] \,,\, P_{H,m} \leq P_{H_\mathrm{out}} \leq P_{H,M}}\; \left\vert \dfrac{\partial f_{m^*}(P_{H_\mathrm{out}})}{\partial P_{H_\mathrm{out}}} \right\vert < +\infty , \\
\zeta & = \inf\limits_{m^* \in [\eta , 1] \,,\, P_{H,m} \leq P_{H_\mathrm{out}} \leq P_{H,M}}\; \alpha_{m^*} \left( P_{H_\mathrm{out}} \right) > 0 .
\end{align*}

We investigate first the case $p=2$. Using Lemmas~\ref{lem 1}-\ref{lem 2}, we obtain for all $t \geq 0$ that:
\begin{equation*}
\dot{V}(t) 
\leq  - 4 \zeta  c_f V(t) + \sqrt{\dfrac{C_f}{c_f}} \left\{
F_1 \vert \dot{m}^*(t) \vert 
+ F_2 \vert \dot{P}_{H_\mathrm{out}}(t) \vert 
\right\} .
\end{equation*}
An integration shows that, for all $t \geq 0$,
\begin{align*}
V(t) 
& \leq e^{-4 \zeta c_f t} V(0) \\
& \phantom{\leq}\, + \dfrac{1}{4 \zeta c_f} \sqrt{\dfrac{C_f}{c_f}} \left\{
F_1 \sup\limits_{\tau \in [0,t]} \left\vert \dot{m}^*(\tau) \right\vert
+ F_2 \sup\limits_{\tau \in [0,t]} \left\vert \dot{P}_{H_\mathrm{out}}(\tau) \right\vert 
\right\} .
\end{align*}
Using the definition of $V$ and the fact that $\sqrt{a+b} \leq \sqrt{a} + \sqrt{b}$ for all $a,b \geq 0$, we obtain that (\ref{eq: ISS P_M_out}) with $p=2$ holds true with $\beta = 2 \zeta c_f > 0$, $C_1 = \sqrt{\dfrac{F_1}{2 \zeta c_f} \sqrt{\dfrac{C_f}{c_f}}} > 0$ and $C_2 = \sqrt{\dfrac{F_2}{2 \zeta c_f} \sqrt{\dfrac{C_f}{c_f}}} > 0$. Now, assuming a cooperative behavior, i.e., $P_{H_\mathrm{out}}(t) \leq P_T(m^*(t))$ for all $t \geq 0$, we obtain from (\ref{eq: def fm*}) that
\begin{align*}
& \vert m(t) - m^*(t) \vert \\
& = \left\vert \dfrac{P_{H_\mathrm{out}}(t)}{P_{M_\mathrm{out}}(t)+P_{H_\mathrm{out}}(t)} - \dfrac{P_{H_\mathrm{out}}(t)}{\sqrt{f_{m^*(t)} \left( P_{H_\mathrm{out}}(t) \right)}+P_{H_\mathrm{out}}(t)} \right\vert \\
& = \vert P_{H_\mathrm{out}}(t) \vert \dfrac{\left\vert P_{M_\mathrm{out}}(t) - \sqrt{f_{m^*(t)} \left( P_{H_\mathrm{out}}(t) \right)} \right\vert}{\left\vert P_{M_\mathrm{out}}(t)+P_{H_\mathrm{out}}(t) \right\vert \left\vert \sqrt{f_{m^*(t)} \left( P_{H_\mathrm{out}}(t) \right)}+P_{H_\mathrm{out}}(t) \right\vert} \\
& \leq \dfrac{P_{H,M}}{( \sqrt{c_f} + P_{H,m} )^2} \left\vert P_{M_\mathrm{out}}(t) - \sqrt{f_{m^*(t)} \left( P_{H_\mathrm{out}}(t) \right)} \right\vert .
\end{align*}
Using (\ref{eq: ISS P_M_out}), it is sufficient to note based on the latter identity applied at $t=0$ that
\begin{align*}
& \left\vert P_{M_\mathrm{out},0} - \sqrt{f_{m^*(0)} \left( P_{H_\mathrm{out}}(0) \right)} \right\vert \\
& = \dfrac{\left\vert P_{M_\mathrm{out},0}+P_{H_\mathrm{out}}(0) \right\vert 
\left\vert \sqrt{f_{m^*(0)} \left( P_{H_\mathrm{out}}(0) \right)}+P_{H_\mathrm{out}}(0) \right\vert}{\vert P_{H_\mathrm{out}}(0) \vert} \\
& \phantom{=}\; \times  \vert m(0) - m^*(0) \vert \\
& \leq \dfrac{\left(\sqrt{C_f}+P_{H,M}\right)^2}{P_{H,m}} \vert m(0) - m^*(0) \vert ,
\end{align*}
to conclude the existence of constants $C_3,C_4,C_5>0$, that are independent of the initial condition and the signals $P_{H_\mathrm{out}},m^*$, such that (\ref{eq: ISS m}) with $p=2$ holds true.

We now consider the case $p=1$. From (\ref{eq: dot_V}) and using the Young's inequality\footnote{Young's inequality: $ab \leq \dfrac{a^2}{2 \delta} + \dfrac{\delta b^2}{2}$ for all $a,b \geq 0$ and $\delta > 0$.},
\begin{align*}
\dot{V}(t)
& \leq  - 4 \zeta c_f V \\
& \phantom{=}\, + \dfrac{1}{2 \sqrt{c_f}} \left\{
F_1 \vert \dot{m}^* \vert
+ F_2 \vert \dot{P}_{H_\mathrm{out}} \vert
\right\} \left\vert P_{M_\mathrm{out}} - \sqrt{f_{m^*}\left( P_{H_\mathrm{out}} \right)} \right\vert \\
& \leq - 2 \zeta c_f V 
+ \dfrac{1}{8 \zeta c_f^2} \left\{
F_1^2 \vert \dot{m}^* \vert^2
+ F_2^2 \vert \dot{P}_{H_\mathrm{out}} \vert^2
\right\} .
\end{align*}
The remaining of the proof for the case $p=1$ is similar to the case $p=2$.
\end{pf*}

%%%%%%%%%%%%%%%%%%%%%%%%%%%%%%%%%%%%%%%%%%%%%%%%%%%%%%%%%%%%%%%%%%%%%%%%%%%%%%%%
\section{Experimental results}\label{sec: results}

The experimental setup involves the e-bike depicted in Fig.~\ref{Ebike} which is mounted securely in a stationary cycling turbo trainer. This allows for rotation of the rear wheel with minimal resistance and, at the same time, to measure the e-bike speed at the wheel.

As discussed in Section~\ref{sec: control design}, the amount of power the motor delivers is controlled by means of the control input  $Y(t)$. But due to the physical characteristics of the motor, the amount of motor power delivered for a given value of $Y(t)$ depends on the human power input (discussed further in \citep{sweeney2018context}). As such, we introduce a secondary control algorithm to get the actual motor power to converge to the target motor power. This secondary control algorithm defined a regulation error $e_k = P_{M_\mathrm{out - Target}} - P_{M_\mathrm{out - Actual}}$ which was used to update the control input $Y_{k+1} = Y_{k} - \gamma e_k$, where $\gamma > 0$ is a proportional gain. The sampling period of the measurement system on the bike was limited to 1 second which limited the sampling period of the secondary control algorithm. If the sampling period could be faster, it is reasonable to expect the convergence of the actual motor power to the target motor power would happen faster. 

The aim of the experiments reported in this section is the validation of the control law described in the previous section. The experiments test the control strategy for a human operating in the both competitive and cooperative regions of the characteristic depicted in Fig.~\ref{fig: function f coop vs comp}. The experiments were completed for a range of values of $m^*$, emulating different pollution levels. Indeed, as $m^*$ stands for the desired value of the power ratio (\ref{eq: m - naive approach}), a low value of $m^*$ corresponds to a e-bike operating mostly in electrical mode. This is particularly suited for highly polluted areas to reduce the effort of the cyclist, hence reducing its ventilation rate. Conversely, a high value of $m^*$ corresponds to a e-bike operating mostly with the power delivered by the cyclist, which is intended for lowly polluted areas. The policy of selecting $P_T(m^*)$ was that for higher pollution areas (i.e. lower values of $m^*$), $P_T$ should also be lower to try to keep the cyclist's effort lower to avoid breathing in high levels of pollution. Other policies are possible for selecting this threshold depending on the desired outcome.

\subsection{Experimental validation of controller design}

\subsubsection{Disturbance rejection for a constant reference $m^*$}

\begin{figure}[htb]
\centering
\includegraphics[width=3.3in]{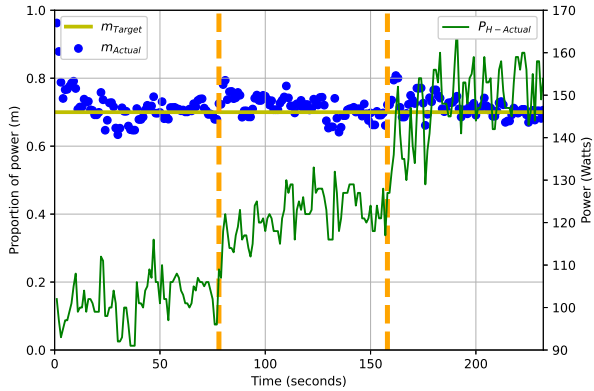}
\caption{Disturbance rejection w.r.t. variations of the power delivered by the cyclist for $m^* = 0.7$. The vertical lines denote a change of trend regarding the effort provided by the cyclist.}
\label{fig:cyclistDisturbanceRejection}
\end{figure}

First, the system has been tested with a constant reference input $m^*$ to evaluate the tracking capability of the closed-loop system in the presence of an input disturbance, namely, a sudden change in the power delivered by the human $P_{H_\mathrm{in}}$. We consider here the case of a cooperative scenario, i.e., the power delivered by the human remains below the threshold $P_T(m^*)$. The experimental result is depicted in Fig.~\ref{fig:cyclistDisturbanceRejection}.

It is observed that the constant reference $m^* = 0.7$ can be tracked satisfactorily in the presence of small variations of the power $P_{H_\mathrm{in}}(t)$ delivered by the human. The tracking is also ensured in the presence of a larger change in average human power from approximately $100\,\mathrm{W}$ to $150\,\mathrm{W}$, which occurs from time $t = 80\,\mathrm{s}$ to time $t = 160\,\mathrm{s}$. Indeed, a sudden increase of the power delivered by the human induces first a significant increase of the ratio $m$. Nevertheless, the control algorithm successfully increases the power delivered by the motor to regulate the value of the ratio $m$ to track the constant reference value $m^* = 0.7$.

\begin{figure}[htb]
\centering
\subfigure[Tracking performance and power inputs. The vertical lines denote a change of trend regarding the effort provided by the cyclist.]{
\includegraphics[width=3.3in]{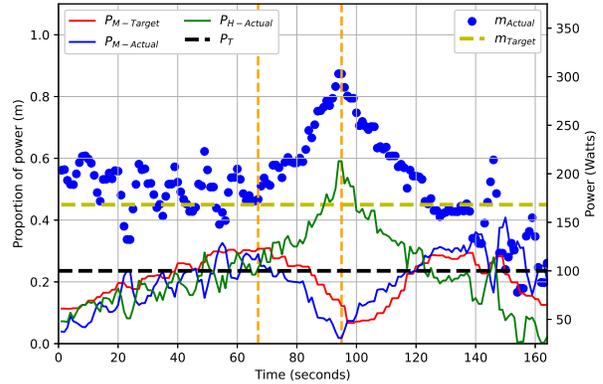}
}
\subfigure[Motor power versus human power]{
\includegraphics[width=3.3in]{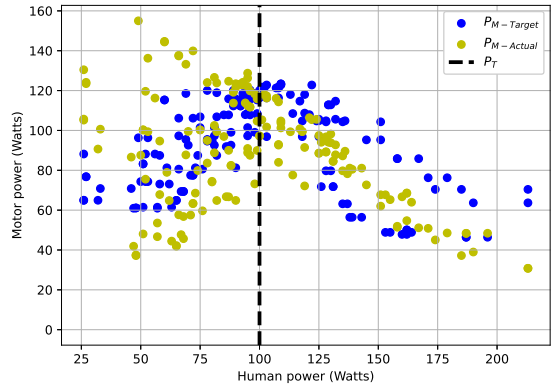}
}
\caption{Reference tracking for $m^* = 0.45$}
\label{fig:cyclistRegZero45}
\end{figure}

\begin{figure}[htb]
\centering
\subfigure[Tracking performance and power inputs. The vertical lines denote a change of trend regarding the effort provided by the cyclist.]{
\includegraphics[width=3.3in]{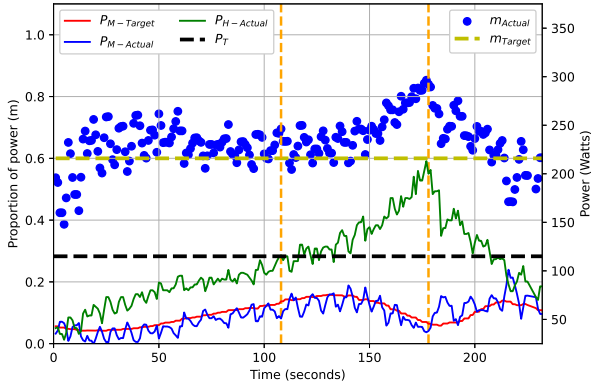}
}
\subfigure[Motor power versus human power]{
\includegraphics[width=3.3in]{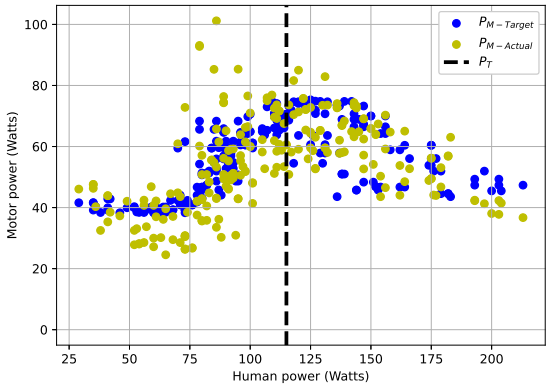}
}
\caption{Reference tracking for $m^* = 0.6$}
\label{fig:cyclistRegZero6}
\end{figure}

We now illustrate with some experiments the main feature of our control strategy, namely: the management of cooperative and competitive modes. The results presented are for $m^*=0.45$ and $m^*=0.6$ in Fig.~\ref{fig:cyclistRegZero45} and Fig.~\ref{fig:cyclistRegZero6}, respectively. It can be clearly seen in the first part of both sets of results that $m$ is regulated to the target value $m^*$ in the cooperative region, but when the cyclist increases their effort beyond the threshold $P_T(m^*)$ and into the competitive region, the target motor power starts to be gradually reduced as the human continues to increase their effort further into the competitive region. As as expected in the competitive region, $m$ is no longer regulated to $m^*$. When the human starts to reduce their effort, the target motor power starts to be increased again.

Experiments were also completed for higher values of $m^*$ such as $m^*=0.9$, but in this situation the human should be providing almost all of the power anyway so there is little difference trying to regulate the motor power. Additionally, experiments were completed for values of $m^*$ lower than 0.45 but the system became unstable at this point due to the large sampling period which was limited by the measurement system on the bike.

\subsubsection{Tracking of a time-varying reference signal $m^*$}

We now report the test of the tracking performance with respect to a time-varying $m^*$. The reference signal $m^*$ is initially set as a constant value of $0.75$. Then, it is smoothly transitioned over the range of time from $t = 70\,\mathrm{s}$ to $t = 100\,\mathrm{s}$ to the constant final value of $0.45$.

\begin{figure}[htb]
\centering
\includegraphics[width=3.3in]{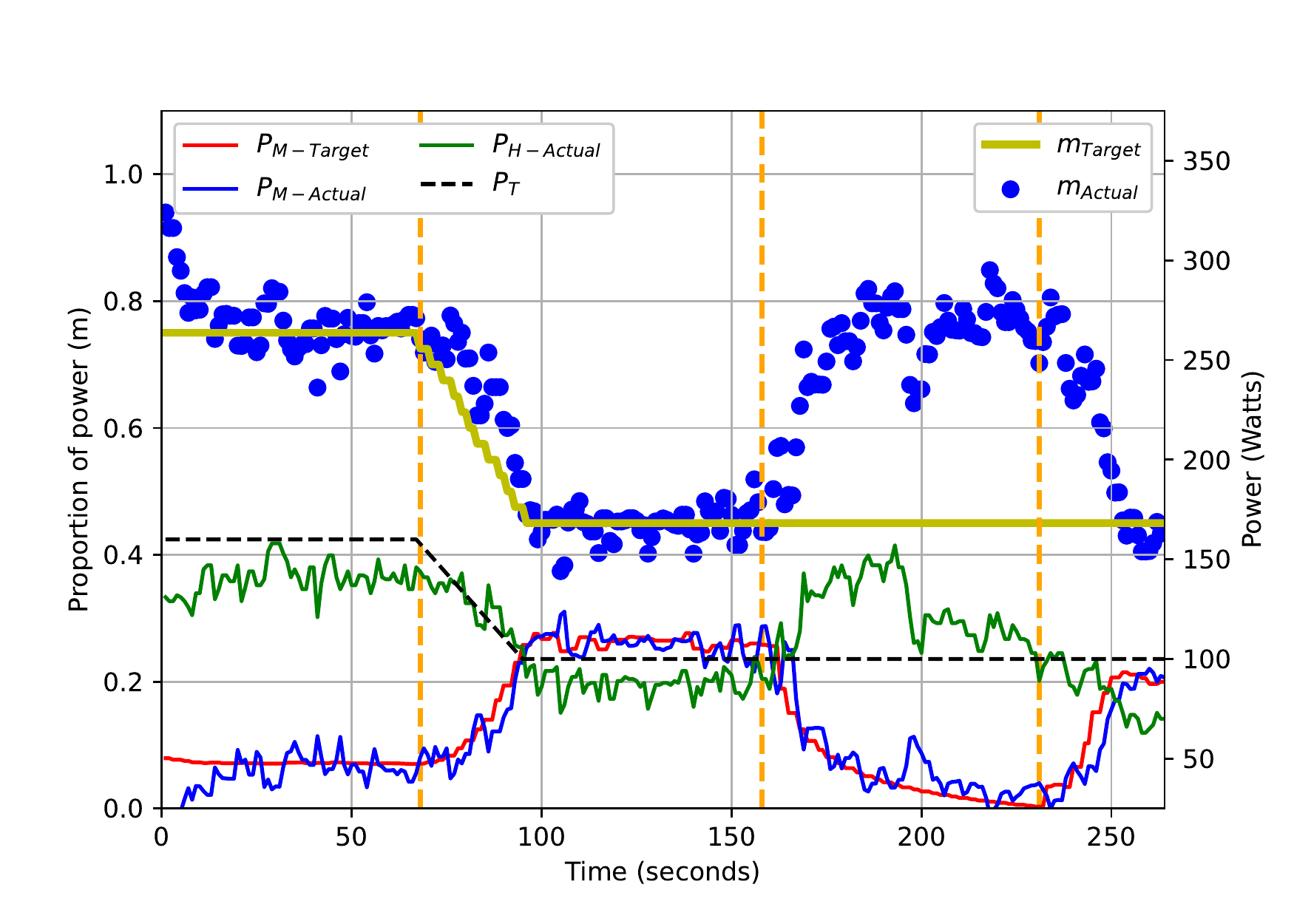}
\caption{Reference tracking for a time-varying $m^*$. The vertical lines denote a change of trend regarding the effort provided by the cyclist.}
\label{fig:track_m_star}
\end{figure}

The experimental results are depicted in Fig.~\ref{fig:track_m_star}. Over the range of time $t = 0\,\mathrm{s}$ to $t = 160\,\mathrm{s}$, the cyclist maintains their effort (essentially) below the threshold $P_T(m^*)$. Hence the power ration $m$ is regulated in a satisfactory manner to the time-varying reference $m^*$. From time $t = 160\,\mathrm{s}$, the cyclist enters into the competitive mode, inducing a significant reduction of the power delivered by the motor. The cooperative mode is resumed from time $t = 230\,\mathrm{s}$ until the end of the experiment.

\subsection{Indirect control of the ventilation rate of the cyclist}

In this test case, we show that the control strategy investigated in this paper can be successfully used to perform an indirect control of the breathing rate of the cyclist. In this experiment, the ventilation rate (VR) of the cyclist is monitored by the spirometry equipment \textit{COSMED Spiropalm 6MWT} with a \textit{COSMED} VO\textsubscript{2} max digital flowmeter. In this setting, the reference signal $m^*$ is generated depending on a time-varying artificial environmental pollution signal.   Specifically, large values of $m^*$ (i.e. close to $1$) correspond to low pollution levels. In this case most of the power is delivered by the human which can freely adjust his effort. Conversely, small values of $m^*$ (i.e. close to the lower bound $\eta > 0$) correspond to high pollution levels. In this case, the motor assists the cyclist in order to limit the effort provided by the human.

\begin{figure}[htb]
\centering
\subfigure[Tracking performance and ventilation rate]{
\includegraphics[width=3.3in]{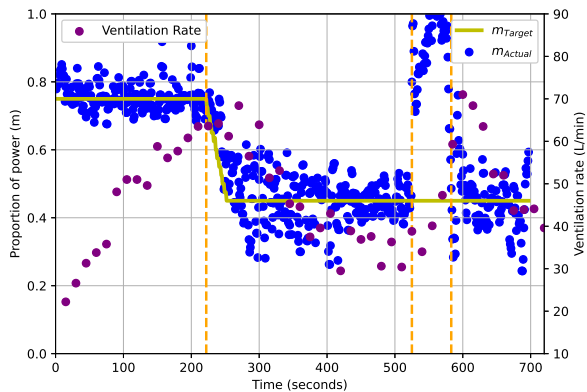}
\label{fig: step_changes_ventilation - ventilation rate}
}
\subfigure[Tracking performance and power]{
\includegraphics[width=3.3in]{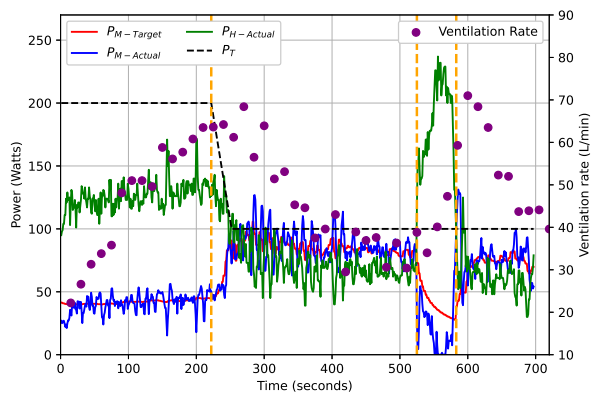}
\label{fig: step_changes_ventilation - power inputs}
}
\caption{Tracking a time-varying $m^*$ which is set according to an artificial pollution level signal, with effect on cyclist's ventilation rate. The vertical lines denote a change of trend regarding the effort provided by the cyclist.}
\label{fig:step_changes_ventilation}
\end{figure}

The experimental results are depicted in Fig.~\ref{fig:step_changes_ventilation}. We observe that the reference signal $m^*$ is adequately tracked by the actual ratio $m$ as long as the cyclist remains in cooperative mode (time $t \leq 520\,\mathrm{s}$). During the first part of the experiment (time $t \leq 520\,\mathrm{s}$), most of the poser is provided by the cyclist with a ratio $m$ around $0.75$. This progressively induces an increase of the ventilation rate of the cyclist. Around time $t = 220\,\mathrm{s}$, the reference value $m^*$ of the ratio is lowered to the value $m^* = 0.45$, indicating a relative increase of the power provided by the motor compared to the one delivered by the cyclist. This gradually reduces the effort of the cyclist, hence reduces his ventilation rate. Such a tendency is stopped in the competitive scenario (times ranging from $t=525\,\mathrm{s}$ to $t=580\,\mathrm{s}$) while restored as soon as the cyclist goes back to a cooperative mode (for times larger than $t=580\,\mathrm{s}$).

%%%%%%%%%%%%%%%%%%%%%%%%%%%%%%%%%%%%%%%%%%%%%%%%%%%%%%%%%%%%%%%%%%%%%%%%%%%%%%%%
\section{Conclusion}\label{sec: conclusion}

We presented a human-in-the-loop control system for deployment on an e-bike.  Specifically, we introduced a control algorithm that is able not only to manage route disturbances but also to manage the interactions with the cyclist.  The basis of the control was a pitchfork bifurcation system, augmented using ideas from gain-scheduling.  After having investigated the robustness of the proposed approach, we experimentally evaluated the effectiveness of our system.

%%%%%%%%%%%%%%%%%%%%%%%%%%%%%%%%%%%%%%%%%%%%%%%%%%%%%%%%%%%%%%%%%%%%%%%%%%%%%%%%
%\bibliographystyle{elsarticle-harv}   
%\bibliography{mybibfile}

\end{document}